\documentclass[twocolumn]{autart}    

\usepackage{graphicx}      
\usepackage[T1]{fontenc}
\usepackage[utf8]{inputenc}
\usepackage{float}
\usepackage{amssymb}
\usepackage{amsmath,color}
\usepackage{subcaption}
\usepackage{stmaryrd}
\usepackage{mathbbol}
\usepackage{enumerate}
\usepackage{mathtools}
\usepackage{paralist}
\usepackage{bm}
\usepackage{tikz}
\usepackage[normalem]{ulem} 
\usepackage{etoolbox}

\newcommand*\MyScale{1}

\tikzset{every picture/.style={scale=\MyScale,}}

\usepackage{pgfplots}
\pgfplotsset{compat=newest}
\pgfplotsset{plot coordinates/math parser=false}

\newlength\fheight 
\newlength\fwidth  
\newcommand{\R}{\mathbb{R}}

\newcommand\numeq[2]%
  {\stackrel{\scriptscriptstyle(\mkern-1.5mu#1\mkern-1.5mu)}{#2}} 

\newtoggle{proofAutomatica}
\toggletrue{proofAutomatica}

\begin{document}

\begin{frontmatter}

\title{Exponential convergence of recursive least squares with forgetting factor for multiple-output systems} 

\author[UCSD]{Sven Br\"{u}ggemann}\ead{sbruegge@eng.ucsd.edu},    
\author[UCSD]{Robert R. Bitmead}\ead{rbitmead@eng.ucsd.edu}

\address[UCSD]{Mechanical \&\ Aerospace Engineering Department, University of California, San Diego, CA 92093-0411, USA}


\begin{abstract}                          
We extend results of the recursive-least-squares-with-forgetting-factor identifier for single-input-single-output systems to the multiple-output case by deriving the corresponding minimized objective function and by showing exponential convergence of the estimation error.
\end{abstract}

\end{frontmatter}

\section{Introduction}
This work revolves around the popular recursive-least-squares-with-forgetting-factor (RLSFF) estimation algorithm for multiple-output (MO) systems. The idea of a recursive formulation goes back to \cite{goodwin1977}, where it is shown that the estimate converges to a value which minimizes the sum of output errors, whereby the outputs' weights reduce exponentially with respect to their time of measurement; hence, the name forgetting factor. These results are established for single-input-single-output (SISO) systems and extended to MO systems in \cite{bernstein2019}. Since, the RLSFF algorithm has been proposed in many different variations, particulalry for the SISO case (see e.g. \cite{Dasgupta}, \cite{SALGADO}, \cite{Leung2005}). Also, exponential convergence of the estimation error under the influence of the RLSFF identifier from \cite{goodwin1977} is shown in \cite{JOHNSTONE1982}. The latter work contains a deterministic analysis considering SISO systems without noise.

Despite the popularity of the RLSFF algorithm, to our best knowledge, a concise convergence result of the RLSFF for the MO case has not been published. Given the literature on adaptive control algorithms using a recursive-least-squares identifier in one way or the other, whose analysis is restricted to SISO systems (see e.g. \cite{Astrom1994},  \cite{Leung1999}, \cite{Paleologu2014}, \cite{Vahidi2005}), especially within the framework of model predictive control (see e.g. \cite{HEIRUNG201564}, \cite{HEIRUNG2019128}, \cite{HEIRUNG2017340}, \cite{marafioti2014}); and the recent article \cite{bernstein2019} dealing with MO systems, we are convinced that SISO-equivalent results for the MO case are of interest and facilitate further research in adaptive control; especially in the area of adaptive model predictive control, which has attracted relatively little attention \cite{MAYNE20142967}. 

Our note therefore closes the gap between the scalar and the MO case by extending the results from \cite{bernstein2019} and \cite{JOHNSTONE1982}. Under the usual assumption of persistence of excitation, we show that the RLSFF algorithm for MO systems minimizes a criterion for which the multiple outputs are weighted not only by their time of measurement (through the forgetting factor), but also among each other (via a user-defined weighting matrix). Furthermore, we demonstrate that the estimation error convergences exponentially to the origin.

\section{Main results}
Consider the system
\begin{align}\label{eq:y_k}
y_{k+1}=\psi^T_k\theta,
\end{align}
where $ y_k\in\R^n$ is the measurable output, $\psi_k\in\R^{m\times n}$ the regressor matrix, respectively at time $k$, and $\theta\in\R^m$ the constant parameter vector. Suppose that the parameter vector is uncertain and, thus, to be estimated. Toward this end, consider the following RLSFF algorithm.
\begin{align}\label{eq:thetaHatRecursion}
\hat \theta_{k+1}=\hat \theta_{k}+P_{k-1}\psi_{k}D_k^{-1}\left( y_{k+1}-\psi_{k}^T\hat \theta_{k}\right),
\end{align}
where $D_k=\lambda T+\psi_{k}^TP_{k-1}\psi_{k}$ with $T=T^T>0\in\R^{n\times n}$ and constant forgetting factor $\lambda\in(0,1)$, and
\begin{align}\label{eq:Pk_recursion}
P_{k+1}&=\lambda^{-1}\Big(I-P_{k}\psi_{k+1}D_{k+1}^{-1}\psi_{k+1}^T\Big)P_{k},
\end{align}
or equivalently,
\begin{align}\label{eq:PinvRecursion}
P_{k-1}^{-1}&=\lambda P^{-1}_{k-2}+\psi_{k-1}T^{-1}\psi_{k-1}^T.
\end{align}
where $P_{-1}\in\R^{m\times m}$ is symmetric positive definite. Define persistence of excitation as follows.
\begin{defn} \label{def:PE}
The matrix sequence $\{\psi_k\}$ is said to be persistently exciting (PE) if for some constant $S$ and all $j$ there exist positive constants $\alpha$ and $\beta$ such that
\begin{align*}
0<\alpha I \leq \sum_{i=j}^{j+S}\psi_i\psi_i^T\leq\beta I<\infty.
\end{align*}
\end{defn}

The following lemma is akin to the result in \cite{bernstein2019}.
\begin{lem}\label{lem:minCostFun}
Suppose $\{\psi_k\}$ is PE. Then, for $k\geq S$, the algorithm in \eqref{eq:thetaHatRecursion} and \eqref{eq:Pk_recursion} converges to the value $\theta$ which minimizes
\begin{align*}
\sum_{i=1}^k \lambda^{k-i}|y_i-\psi^T_{i-1} \theta|^2_{T^{-1}}.
\end{align*}.
\end{lem}
\begin{pf}
\iftoggle{proofAutomatica}
{
The proof is analogous to that of \cite[Theorem 2]{bernstein2019} and hence omitted for brevity.
}
{
Let
\begin{align*}
L(\theta)&=\frac{1}{2} \sum_{i=1}^k\lambda^{k-i}|y_i-\psi_{i-1}^T\theta|^{2}_{T^{-1}}\\
&=\frac{1}{2} \sum_{i=1}^k \lambda^{k-i}\left(y_i^TT^{-1}y_i-2\theta^T\psi_{i-1}T^{-1}y_i\right.\\
& \quad \quad \quad \quad \quad \quad \left. +\theta^T\psi_{i-1}T^{-1}\psi_{i-1}^T\theta\right).
\end{align*}
Differentiating this yields
\begin{align*}
\frac{\partial L(\theta)}{\partial \theta}&=\sum_{i=1}^k \lambda^{k-i}\left(\psi_{i-1}T^{-1}\psi_{i-1}^T\theta-\psi_{i-1}T^{-1}y_i\right)\\
&=\sum_{i=1}^k \lambda^{k-i}\psi_{i-1}T^{-1}\left(\psi_{i-1}^T\theta-y_i\right)\\
&=\left(\sum_{i=1}^k \lambda^{k-i}\psi_{i-1}T^{-1}\psi_{i-1}^T\right)\theta\\
&\quad -\sum_{i=1}^k \lambda^{k-i}\psi_{i-1}T^{-1}y_i.
\end{align*}
Denote $\hat \theta_k$ as the value of $\theta$ which satisfies the equation above set to zero. Then, for $k\geq S$, one obtains
\begin{align}\label{eq:theta minimizer}
\hat \theta_k=\underbrace{\left(\sum_{i=1}^k \lambda^{k-i}\psi_{i-1}T^{-1}\psi_{i-1}^T\right)^{-1}}_{\coloneqq P_{k-1}}\sum_{i=1}^k \lambda^{k-i}\psi_{i-1}T^{-1}y_i,
\end{align}
where the matrix $P_{k-1}^{-1}$ can be written in a recursive manner:
\begin{align}\label{eq:PinvRecursion}
P_{k-1}^{-1}&=\sum_{i=1}^k \lambda^{k-i}\psi_{i-1}T^{-1}\psi_{i-1}^T\nonumber\\
&=\lambda P^{-1}_{k-2}+\psi_{k-1}T^{-1}\psi_{k-1}^T.
\end{align}
Using the Woodbury matrix identity leads to the following recursive expression:
\begin{align}\label{eq:Precursion}
P_{k-1}&=\lambda^{-1}P_{k-2}\nonumber\\
&\quad -\lambda^{-1}P_{k-2}\psi_{k-1}\left(T+\psi_{k-1}^T\lambda^{-1}P_{k-2}\psi_{k-1}\right)^{-1}\nonumber\\
&\quad \psi^T_{k-1}\lambda^{-1}P_{k-2}\nonumber\\
&=\lambda^{-1}\left(I \vphantom{\left(T^T \right)^{-1}} \right.\nonumber\\
&\quad\left.-\lambda^{-1} P_{k-2}\psi_{k-1}\left(T+\psi_{k-1}^T\lambda^{-1}P_{k-2}\psi_{k-1}\right)^{-1}\psi_{k-1}^T\right)\nonumber\\
&\quad P_{k-2}\nonumber\\
&=\lambda^{-1}\Big(I \nonumber\\
&\quad-P_{k-2}\psi_{k-1}\underbrace{\left(\lambda T+\psi_{k-1}^TP_{k-2}\psi_{k-1}\right)^{-1}}_{\coloneqq D_{k-1}^{-1}}\psi_{k-1}^T\Big)P_{k-2},
\end{align}
which is equivalent to \eqref{eq:Pk_recursion}.

Toward the goal of a recursive expression for the parameter vector $\hat \theta_k$, consider
\begin{align*}
\hat \theta_k&\numeq{\ref{eq:theta minimizer}}{=}P_{k-1}\sum_{i=1}^k \lambda^{k-i}\psi_{i-1}T^{-1}y_i\nonumber\\
&=P_{k-1}\left(\sum_{i=1}^{k-1} \lambda^{k-i}\psi_{i-1}T^{-1}y_i + \psi_{k-1}T^{-1}y_k\right)\nonumber\\
&\numeq{\ref{eq:theta minimizer}}{=} P_{k-1}\left(\lambda P_{k-2}^{-1}\hat \theta_{k-1}+ \psi_{k-1}T^{-1}y_k\right)\nonumber\\
&\numeq{\ref{eq:Precursion}}{=}\hat \theta_{k-1}+P_{k-2}\psi_{k-1}D_{k-1}^{-1}\left(y_k-\psi_{k-1}^T\hat \theta_{k-1}\right),
\end{align*}
which derives \eqref{eq:thetaHatRecursion}.
}
\hfill$\blacksquare$
\end{pf}
We have thus presented the objective which is minimized by the RLSFF algorithm in \eqref{eq:thetaHatRecursion} and \eqref{eq:Pk_recursion}. In contrast to the existing literature, we explicitly incorporate a weight for each respective output via the matrix $T$. We now wish to obtain an exponentially stable estimation error
\begin{align*}
\tilde \theta_{k}&=\theta-\hat \theta_{k}.
\end{align*}
The corresponding theorem as an extension of \cite{JOHNSTONE1982} follows.

\begin{thm}\label{lem:tilde theta convergence}
Suppose $\{\psi_k\}$ is PE. Then, for any initial condition $\tilde \theta_0$, the estimation error $\tilde \theta_k$ converges exponentially to $\theta$, i.e. for any $\tilde \theta_0$ there exist $\gamma>0$ such that for all $k\geq S$
\begin{align*}
|\tilde \theta_k|^2\leq \gamma \lambda^{k} |\tilde \theta_0|^2.
\end{align*}
\end{thm}
\begin{pf}
The proof is divided into two parts. Part I establishes a lower bound on $P^{-1}_k$. This is then used in part II to show exponential stability of the estimation error.

\textit{Part I}:
Recollect that if $B$ is symmetric, then for any matrix $A$,
\begin{align*}
ABA^T\geq \lambda_{min}(B)AA^T.
\end{align*}
This follows by definition of a positive definite matrix;
\begin{align*}
x^T(ABA^T&- \lambda_{min}(B)AA^T)x\\
&=x^TABA^Tx- \lambda_{min}(B) x^TAA^Tx\\
&=(A^Tx)B(A^Tx)- \lambda_{min}(B) x^TAA^Tx\\
&\geq \lambda_{min}(B) (A^Tx)A^Tx -  \lambda_{min}(B) x^TAA^Tx\\
&=0,
\end{align*}
where we use the fact that $\lambda_{min}(B)|x|^2\leq |x|_B^2$. It follows that if $\{\psi_k\}$ is PE, then
\begin{align*}
P^{-1}_{j-1}+\dots+P^{-1}_{j+S-1}&\numeq{\ref{eq:PinvRecursion}}{\geq} \sum_{k=j}^{j+S}\psi_{k-1}T^{-1}\psi_{k-1}^T\\
&~\geq \lambda_{min}(T^{-1})\alpha I.
\end{align*}
for all $k\geq S$. Following \cite[Lemma 1]{JOHNSTONE1982} leads to the lower bound
\begin{align}\label{eq:PinvLB}
P^{-1}_{k-1}\geq \frac{\lambda_{min}(T^{-1})\alpha(\lambda^{-1}-1)}{\lambda^{-(S+1)}-1}>0
\end{align}
for all $k\geq S$.

\textit{Part II}:
By \eqref{eq:thetaHatRecursion}, one can write recursively the estimation error as
\begin{align}\label{eq:theta tilde recursion}
\tilde \theta_{k+1}=\left(I-P_{k-1}\psi_{k}D_k^{-1}\psi^T_{k}\right)\tilde \theta_{k}.
\end{align}
Consider the Lyapunov function candidate
\begin{align}\label{eq:Wk}
W_k=\tilde \theta_k^T P_{k-1}^{-1}\tilde \theta_k.
\end{align}
Then, using the recursions in \eqref{eq:PinvRecursion} and \eqref{eq:theta tilde recursion} yields
\begin{align}\label{eq:Wk+1-Wk}
W_{k+1}-W_k&=\tilde \theta_{k+1}^T P_{k}^{-1}\tilde \theta_{k+1}-\tilde \theta_k^T P_{k-1}^{-1}\tilde \theta_k\nonumber \\
&=\tilde \theta_k^T\left[(\lambda-1)P_{k-1}^{-1}-\lambda\psi_k D_k^{-1}\psi_k^T+C\right]\tilde \theta_k,
\end{align}
where
\begin{align*}
C&=\psi_k\left[T^{-1}-D_k^{-1}\psi_k^TP_{k-1}\psi_kT^{-1} -\lambda D_k^{-1}\right.\nonumber\\
&\quad -T^{-1}\psi_k^T P_{k-1}\psi_kD_k^{-1}+\lambda D_k^{-1}\psi_k^TP_{k-1}\psi_kD_k^{-1}\nonumber \\
&\quad \left.+D^{-1}_k\psi^T_kP_{k-1}\psi_kT^{-1}\psi_k^T P_{k-1}\psi_kD_k^{-1}\right]\psi_k^T.
\end{align*}
We demonstrate now that $C$ is equal to the zero matrix. To this end, multiply the inner term of $C$ by $T^{1/2}$ from both sides respectively and obtain
\begin{align}\label{eq:proof T C T}
I&-T^{1/2}D_k^{-1}\psi_k^TP_{k-1}\psi_kT^{-1/2}-\lambda T^{1/2}D_k^{-1}T^{1/2} \nonumber\\
&-T^{-1/2}\psi_k^TP_{k-1}\psi_kD_k^{-1}T^{1/2} \nonumber\\
&+\lambda T^{1/2}D_k^{-1}\psi_k^TP_{k-1}\psi_kD_k^{-1}T^{1/2} \nonumber\\
&+T^{1/2}D_k^{-1}\psi_k^TP_{k-1}\psi_kT^{-1}\psi_k^TP_{k-1}\psi_k D_k^{-1}T^{1/2}.
\end{align}
By defining
\begin{align}\label{eq:Dbar}
\bar \psi&=\psi_kT^{-1/2}\nonumber\\
\bar D&=T^{-1/2}D_kT^{-1/2}=\lambda I+\bar\psi^TP_{k-1}\bar \psi
\end{align}
one has that
\begin{align*}
T^{1/2}D_k^{-1}&=(D_kT^{-1/2})^{-1}=(\lambda T^{1/2}+\psi_k^TP_{k-1}\psi_kT^{-1/2})^{-1}\\
&=\left(T^{1/2}(\lambda+T^{-1/2}\psi_k^TP_{k-1}\psi_kT^{-1/2})\right)^{-1}\\
&=\left(T^{1/2}(\lambda+\bar \psi^TP_{k-1}\bar \psi)\right)^{-1}\\
&=\bar D^{-1}T^{-1/2},
\end{align*}
so that for \eqref{eq:proof T C T} it follows:
\begin{align*}
I&-\bar D^{-1}\bar \psi^TP_{k-1}\bar \psi-\lambda\bar D^{-1}\\
&-\bar \psi^TP_{k-1}\bar \psi\bar D^{-1}+\lambda \bar D^{-1}\bar \psi^TP_{k-1}\bar \psi\bar D^{-1}\\
&+\bar D^{-1}\bar \psi^TP_{k-1}\bar \psi\bar \psi^TP_{k-1}\bar \psi\bar D^{-1}.
\end{align*}
Observe that this can be reformulated as
\begin{align*}
I&-\bar D^{-1}\overbrace{\left(\bar \psi^TP_{k-1}\bar\psi+\lambda I\right)}^{\numeq{\ref{eq:Dbar}}{=}\bar D}\\
&\quad +\bar D^{-1}\left(-\bar D\bar \psi^TP_{k-1}\bar \psi+\lambda\bar \psi^TP_{k-1}\bar \psi\right.\\
& \quad \quad \quad \quad \quad  \left.+\bar\psi^TP_{k-1}\bar\psi\bar\psi^TP_{k-1}\bar\psi\right)\bar D^{-1},
\end{align*}
which is clearly zero and, thus, so is $C$.

Therefore, the difference related to the Lyapunov function candidate in \eqref{eq:Wk+1-Wk}
\begin{align*}
W_{k+1}-W_k&=\tilde \theta_k^T\left[(\lambda-1)P_{k-1}^{-1}-\lambda\psi_k D_k^{-1}\psi_k^T\right]\tilde \theta_k\\
&\leq (\lambda-1)\tilde \theta_k^TP_{k-1}^{-1}\tilde \theta_k\\
&=(\lambda-1)W_k,
\end{align*}
so that
\begin{align*}
W_{k+1}\leq\lambda W_k\leq\lambda^{k+1}W_0=\lambda^{k+1}\tilde \theta_0^TP_{-1}^{-1}\tilde \theta_0.
\end{align*}
Finally, combining this inequality with the definition of $W_k$ in \eqref{eq:Wk} and the lower bound of $P_{k-1}^{-1}$ in \eqref{eq:PinvLB} leads to
\begin{align*}
|\tilde \theta_k|^2\leq \underbrace{\frac{\lambda^{-(S+1)}-1}{\lambda_{min}(T^{-1})\alpha(\lambda^{-1}-1)}\lambda_{max}\left(P_{-1}^{-1}\right)}_{\eqqcolon \gamma}\lambda^k|\tilde \theta_k|^2
\end{align*}
for all $k\geq S$.\hfill$\blacksquare$
\end{pf}
\begin{rem} Note that exponential stability of the estimation error implies a bounded error in the presence of bounded additive noise. Thus, for additive bounded noise on system \eqref{eq:y_k}, the estimation error converges to a ball centered on the true parameter vector with a radius proportional to the bound on the disturbance.
\end{rem}
\section{Conclusion}
We have shown that the RLSFF estimation algorithm for MO systems shares the properties of those for single output systems, i.e. it minimizes a similar cost function, where outputs are weighted among one another by a given matrix, and it induces exponential convergence of the estimate to the true parameter vector.

\bibliographystyle{plain}        
\bibliography{autosam}           



\end{document}